\documentclass[a4paper, 12pt]{article}

   \usepackage[intlimits]{amsmath}
   \usepackage{amssymb}
   \usepackage{bbm}
   \usepackage{amsthm}
   \usepackage{mathrsfs}
   \usepackage[nice]{nicefrac}

   \usepackage{color}
   \definecolor{sgreen}{rgb}{0.0862745,0.654902,0.054901}
   \definecolor{sred}{rgb}{0.729411,0.145049,0.145049}

   \def\prt{\partial}
   \def\dx{\,\text dx}
   \def\pP{Problem (P)}
   \def\Pp{Problem (P)}
\def\f{\varphi}

   \def\ee{{\,\text e}}
   \def\dvr{\text{div}}
   \def\no{\{1,\ldots,n\}}

   \def\eps{\varepsilon}

   \def\calV{\mathcal V}

   \numberwithin{equation}{section}

   \newtheorem{defi}{Definition}[section]
   
   \newtheorem{coro}[defi]{Corollary}
   \newtheorem{theo}[defi]{Theorem}
   \newtheorem{lemm}[defi]{Lemma}
   \newtheorem{prop}[defi]{Proposition}

   \newcommand{\heat}[1]{\mathcal S(#1)\,}
   \newcommand{\heatT}[1]{\mathcal T(#1)\,}

\begin{document}

\begin{center}
   {\Large \textbf{Contraction in $L^1$ and large time behavior for a system arising in chemical reactions and molecular motors}}
   \vskip0.2in

      {  Michel\ Chipot
         \footnote{An\-ge\-wand\-te Ma\-the\-ma\-tik, Uni\-ver\-si\-taet Zu\-rich, CH-8057 Zurich\\
            E-mail: {\tt m.m.chipot@math.unizh.ch}},\quad
         Danielle\ Hilhorst
         \footnote{Laboratoire de Math\'ematiques, CNRS and
            Universit\'e de Paris-Sud, 91405 Orsay C\'edex, France;
            E-mail: {\tt Danielle.Hilhorst@math.u-psud.fr}},\\
         David Kinderlehrer
         \footnote{Center for Nonlinear Analysis and Department of Mathematical
            Sciences, Carnegie Mellon University, Pittsburgh, PA 15213;
            E-mail: {\tt davidk@andrew.cmu.edu}},
         Micha{\l} Olech
         \footnote{Instytut Matematyczny Uniwersytetu Wroc{\l}awskiego, pl.
            Grunwaldzki 2/4, 50-384 Wroc{\l}aw, Polska; Laboratoire de
            Math\'ematiques, CNRS Universit\'e de Paris-Sud, 91405 Orsay
            C\'edex, France;
            E-mail: {\tt olech@math.uni.wroc.pl}. \\
            The preparation of this article has been partially supported by a KBN/MNiI grant \texttt{1 P03A 008 30} and a Marie Curie Transfer of Knowledge Fellowship of the European Community's Sixth Framework Programme under the contracts \texttt{MTKD-CT-2004-013389}, DMS 0305794, DMS 0405343, and DMS 0806703. We would also like to mention the support of the Swiss National Science Foundation under the contract \# 20-113287/1.}
}
\end{center}

   \abstract{We prove a contraction in $L^1$ property for the solutions of a nonlinear reaction--diffusion system whose special cases include intercellular transport as well as reversible chemical reactions.  Assuming the existence of stationary solutions we show that the solutions stabilize as $t$ tends to infinity. Moreover, in the special case of linear reaction terms, we prove the existence and the uniqueness (up to a multiplicative constant) of the stationary solution.}



%------------------------------------------------------------------------------%
\section{Introduction}\label{HKOintro}
%------------------------------------------------------------------------------%

We start with two specific reaction-diffusion systems. The first one describes a reversible reaction and the other one a molecular motor. We first consider the reversible chemical reaction (see also Bothe \cite{bothe}, Bothe and Hilhorst \cite{bothe_hilhorst}, Desvillettes and Fellner \cite{desvillettes_fellner} and {\'E}rdi and T{\'o}th \cite{erdi_toth}). It involves a reaction-diffusion system of the form
   \begin{equation}\label{HKOProblemBothe}
   \begin{split}
      u_t=&\,d_1 \Delta u - \alpha k \big(r_A(u)-r_B(v)\big) \quad \text{in} \quad \Omega\times(0,T),\ \ \Omega \subset \mathbb R^d,\\
      v_t=&\,d_2 \Delta v + \beta k \big(r_A(u)-r_B(v)\big) \quad \text{in} \quad \Omega \times \ (0,T),\ \ \Omega \subset \mathbb R^d,
   \end{split}
   \end{equation}
   together with the homogeneous Neumann boundary conditions, where $d_1$,\ $d_2$,\ $\alpha$,\ $\beta$,\ $k$ and\ $T$ are positive constants  and where $\Omega$ is a bounded subset of $\mathbb R^d$ with smooth boundary. Such systems describe, with a suitable choice of the functions $r_A$ and $r_B$, chemical reactions for two mobile species. For example, functions $r_A(u)=u^k,\ r_B(v)=v^m$ correspond to a reversible reaction $kA \rightleftharpoons mB$. Reactions of the type $q_1A_1+\ldots q_kA_k \rightleftharpoons q_1B_1 +\ldots q_mB_m$ can also be described by similar systems with more complicated reactions terms.

   Another model problem is a system  in $d=1$ space dimension and $n$ unknown variables $u_1,\dots,u_n$, $n>1$, for in\-ter\-ce\-llu\-lar tran\-sport, namely
   \begin{equation*}
   \begin{split}
      \frac{\partial u_i}{\partial t} &= \frac{\partial}{\partial x} \left(\sigma \frac{\partial u_i}{\partial x} + u_i \psi_i' \right)\\
      & \qquad \qquad+ \sum_{j=1}^n a_{ij}u_j \quad\text{in}\quad Q_T=[0,1]\times (0,T)\\
      & \sigma \frac{\partial u_i}{\partial x} + u_i \psi_i' = 0 \quad \text{on}\quad \prt Q_T=\{0,1\}\times (0,T),
   \end{split}
   \end{equation*}
   where
    \begin{equation}\label{H0}
    \begin{split}
   &a_{ii} \leq 0,\ a_{ij} \geq  0 \mbox{~~for~~all~~}\ \   i\in\no, i \neq j, \\
   &\sum_{i=1}^n a_{ij}=0 \mbox{~~for~~all~~}  i,j \in\no.
   \end{split}
    \end{equation}
   It models transport via motor proteins in the eukaryotic cell where chemical energy is transduced into directed motion. A derivation of the system from a mass transport viewpoint is given in \cite{chipot_kinderlehrer_kowalczyk}.   For an analysis of the steady state solutions and for further references we refer to \cite{chipot_hastings_kinderlehrer}, \cite{hastings_kinderlehrer_mcleod_1}, \cite{hastings_kinderlehrer_mcleod_2}, and \cite{ps}.

In this paper we study the corresponding system in higher space dimension, namely
\begin{subequations}\label{HKOproblem_con}
   \begin{equation}\label{HKOmain_eq_con}
   \begin{split}
      \frac{\prt u_i}{\prt t}  = \dvr\big(\sigma_i\nabla u_i &+ u_i
         \nabla \psi_i \big) \\
            &+ \alpha_i\bigg( \sum_{j=1}^n \lambda_{ij}
         r_j \big(u_j(x,t),x \big)\bigg) \quad \mbox{~~in~~} Q_T,
   \end{split}
   \end{equation}
where $i\in\no$, and $u_i(x,t):Q_T \to \mathbb R^+$, with
$Q_T=\Omega\times (0,T)$, $\Omega$ an open  bounded subset of
$\mathbb R^d$ with smooth boundary, and $T$ some positive constant. We supplement this
system with the Robin (no-flux) boundary conditions
   \begin{equation}\label{HKOmain_bd_con}
      \sigma_i\frac{\prt u_i}{\prt \nu} + u_i\frac{\prt \psi_i}{\prt \nu} = 0,\quad i\in\no,\quad
      \text{on}\quad \partial \Omega \times (0,T),
   \end{equation}
where $\nu$ is the outward normal vector to $\partial\Omega$, and the initial conditions
   \begin{equation}\label{HKOmain_in_con}
      u_1(x,0) = u_{0,1}(x),\ \ldots\ ,u_{n}(x,0) = u_{0,n}(x),\quad x \in \Omega.
   \end{equation}
\end{subequations}

We assume that the following hypotheses hold
\begin{enumerate}
   \item The constants $\sigma_i$ and $\alpha_i\in\mathbb R$, where $i\in \{1,\ldots,n\}$, are strictly positive;
   \item For $i,j\in\no,\ \lambda_{ii} \leq 0,\ \lambda_{ij} \geq 0$ if $i\neq j$,\ $\sum_{k=1}^n \lambda_{kj}=0$;
   \item for all $i\in\no$, the smooth functions $r_i$ are nondecreasing with respect to\ the first variable; $r_i(0,x)=0$ and we assume that the functions $\psi_i$ are smooth as well;
   \item $u_i(.,0) =  u_{0i} \in C(\overline \Omega),\ u_{0i} \geqslant 0$.
\end{enumerate}

In the linear case of the molecular motors, it amounts to choosing
\begin{equation}\label{HKOrownowaznosc}
   r_{i}(s,x) = s,\ \lambda_{ij} = a_{ij}\ \text{and}\ \alpha_i = 1\ \text{for all}\ i,j \in \no.
\end{equation}

We denote by \Pp~the system \eqref{HKOmain_eq_con} together with the
boundary and initial conditions \eqref{HKOmain_bd_con}, \eqref{HKOmain_in_con}, and admit without proof that \Pp~possesses a unique smooth and bounded solution on each time interval $(0,T]$.
An essential idea for proving the existence of a solution would be to apply the Comparison principle Theorem
\ref{HKOmax_princ_theo} below to deduce that any solution of \Pp~ has to be nonnegative and bounded from above by a
stationary solution.

Finally, we note that because of the boundary conditions
\eqref{HKOmain_bd_con} the quantity
   \begin{equation}\label{HKOniezm}
     \sum_{i=1}^n \frac{1}{\alpha_i} \int_\Omega u_i(x,t) \dx
   \end{equation}
is conserved in time.

The organization of this paper is as follows. In Section 2 we prove a comparison principle for \Pp. The main idea, which permits to show that \Pp~is cooperative, is a change of functions which transforms the Robin boundary conditions into the  homogeneous Neumann boundary conditions. In Section 3 we establish a contraction in $L^1$ property for the corresponding semigroup solution. Let us point out the similarity with an old result due to Crandall and Tartar \cite{crandall_tartar} where they proved in a scalar case that in the presence of a conservation of the integral property such as \eqref{HKOniezm}, a comparison principle such as Theorem \ref{HKOmax_princ_theo} is equivalent to a contraction in $L^1$ property such as the inequality \eqref{HKOWeakContr} below. As far as we know such an abstract result is not known in the case of systems.\\
Section 4 deals with the large time behavior of the solutions. Supposing the existence of a stationary solution, we construct a continuum of stationary solutions and prove that the solutions stabilize as $t$ tends to infinity. Let us mention a result by Perthame \cite{perthame} who proved the stabilization in the case of the two component one-dimensional mo\-le\-cu\-lar motor problem. Finally in Section 5, %we make an extra hypothesis that the constants $a_{ij}$ are such that $a_{ij} \ne 0$ for $i,j \in \no$, and
show the existence and uniqueness (up to a multiplicative constant) of the stationary solution of the molecular motor problem. % in space dimensions $1,\ldots, 4$.

\vskip 10pt
\noindent
{\bf Acknowledgment:} The authors acknowledge the preliminary master thesis work of Aude Brisset about the corresponding two component system. They are grateful to the professors\ Piotr\ Biler, Stuart Hastings, Annick Lesne and Hiroshi\ Matano for very fruitful discussions.

%------------------------------------------------------------------------------%
\section{Comparison principle}
%------------------------------------------------------------------------------%
First, we remark that the system of equations \eqref{HKOmain_eq_con} is cooperative.
However, since nothing is known about the sign of the coefficients
$\displaystyle{\frac{\prt \psi_i}{\prt \nu}}$ in the Robin boundary conditions
(\ref{HKOmain_bd_con}), we cannot decide whether the Problem (P)  %(\ref{HKOproblem_con})
is cooperative. This leads us to perform a change of variables which transforms
the Robin boundary conditions into the homogeneous Neumann boundary conditions.
%------------------------------------------------------------------------------%
\subsection{The change of unknown functions}
%------------------------------------------------------------------------------%

Performing the change of variables
\begin{equation}\label{HKOchange2ev}
   w_i(x,t) = u_i(x,t) \ee^{\nicefrac{\psi_i(x)}{\sigma_i}}, \quad i \in \no,
\end{equation}
we deduce from \eqref{HKOproblem_con} that ${\vec w}:=(w_1, \ldots,w_n)$ satisfies
the parabolic problem
 \begin{equation}\label{HKOmain_eq_Neumann}
   \begin{split}
      &\frac{\prt w_i}{\prt t}  =
      \sigma_i \ee^{\nicefrac{\psi_i(x)}{\sigma_i}}
      \dvr\big(\ee^{\nicefrac{-\psi_i(x)}{\sigma_i}}\nabla w_i \big) \\
            &+ \alpha_i \ee^{\nicefrac{\psi_i(x)}{\sigma_i}}
            \bigg( \sum_{j=1}^n \lambda_{ij}
        r_j \big(w_j(x,t)\ee^{\nicefrac{-\psi_j(x)}{\sigma_j}},x \big)\bigg) \quad \mbox{~~in~~} Q_T,
   \end{split}
   \end{equation}
together with the homogeneous Neumann boundary conditions
\begin{equation}\label{HKOnBdNeumann}
   \frac{\partial w_i}{\partial \nu} = 0, \quad i \in \no,\quad \text{on}\quad \partial \Omega,
\end{equation}
and the initial conditions
\begin{equation}\label{HKOivBdNeumann}
   w_i(x,0) =  u_{0,i}(x) \ee^{\nicefrac{\psi_i(x)}{\sigma_i}}, \quad i \in \no, \quad x \in \Omega.
\end{equation}
In the following, we denote by Problem $P_N$ --- the problem (\ref{HKOmain_eq_Neumann}), (\ref{HKOnBdNeumann}), (\ref{HKOivBdNeumann}). To begin with we define the operators
\begin{equation}\label{HKOCLi}
\begin{split}
   &{\cal L}_i(w_i)= \frac{\prt w_i}{\prt t}  -
      \sigma_i \ee^{\nicefrac{\psi_i(x)}{\sigma_i}}
      \dvr\big(\ee^{\nicefrac{-\psi_i(x)}{\sigma_i}}\nabla w_i \big) \\
            &- \alpha_i \ee^{\nicefrac{\psi_i(x)}{\sigma_i}}
            \bigg( \sum_{j=1}^n \lambda_{ij}
        r_j \big(w_j(x,t)\ee^{\nicefrac{-\psi_j(x)}{\sigma_j}},x \big)\bigg)  \quad \mbox{~~in~~} Q_T.
\end{split}
\end{equation}
We say that $\left(\underline w_1,\ldots,\underline w_n\right)$ is a subsolution of Problem $P_N$ if
\begin{equation}\label{subsolution}
\begin{split}
   {\cal L}_i(\underline w_i)\leqslant & \;0 \quad \text{in}\quad Q_T,\\
   \frac{\prt\underline w_i}{\prt\nu} \leqslant & \;0 \quad \text{on}\quad \prt \Omega \times (0,T),\\
      \underline w_i(x,0)\leqslant & \;w_i(x,0), \quad x \in \Omega
\end{split}
\end{equation}
for all $i\in\no$. We define similarly a supersolution $\left(\overline u_1,\ldots,\overline u_n\right)$ of
Problem $P_N$ by the inequalities
\begin{equation}\label{supersolution}
\begin{split}
  {\cal L}_i(\overline w_i)\geqslant & \;0 \quad\text{in}\quad Q_T,\\
   \frac{\prt\overline w_i}{\prt\nu} \geqslant & \;0\quad\text{on}\quad \partial\Omega\times (0,T),\\
   \overline w_i(x,0)\geqslant & \;w_i(x,0),\quad x \in \Omega.
\end{split}
\end{equation}
The following comparison theorem holds (\cite{alikakos_hess_matano}, \cite{protter_weinberger}).
\begin{theo}\label{comptheo}
   Let $(\underline w_1,\ldots,\underline w_n)$ and $(\overline w_1,\ldots,\overline w_n)$, be a sub\,- and a su\-per\,-\,so\-lution, respectively, for the operators ${\cal L}_j$  defined by \eqref{HKOCLi} with $j\in\no$, which means that \eqref{subsolution} and \eqref{supersolution} hold for $i\in\no$. Then $\underline w_i\leqslant \overline w_i$ in $Q_T$. Moreover, for all $i \in \no$ such that $\underline w_i \leqslant \overline w_i$ and $\underline w_i \not \equiv \overline w_i$ on $\{t=0\}\times\Omega$ then $\underline w_i < \overline w_i$ in $Q_T$. \hfill $\blacksquare$
\end{theo}\medskip
This comparison theorem immediately translates into a comparison theorem for solutions of the original \Pp.
For all $i \in \no$, we define the operators
\begin{equation}\label{HKOLi}
\begin{split}
   L_i(u_i)=(u_i)_t-\dvr\big(&\sigma_i\nabla u_i + u_i\nabla \psi_i \big)\\
   &-\alpha_i \bigg(\sum_{j=1}^n \lambda_{ij}\,r_j\,(u_j,x)\bigg)\quad\text{in}\quad Q_T.
\end{split}
\end{equation}
The following result holds.
\begin{theo}\label{HKOmax_princ_theo}
   Let $(\underline u_1,\ldots,\underline u_n)$ and $(\overline u_1,\ldots,\overline u_n)$, be a sub\,- and a su\-per\,-\,so\-lution,  respectively, for the operators $L_j$, defined by \eqref{HKOLi} with $j\in\no$. Then $\underline u_i\leqslant \overline u_i$ in $Q_T$. Moreover, for all $i \in \no$ such that $\underline u_i \leqslant \overline u_i$ and $\underline u_i \not \equiv \overline u_i$ on $\{t=0\}\times\Omega$ then $\underline u_i < \overline u_i$ in $Q_T$. \hskip 5pt $\blacksquare$
\end{theo}\medskip

Next we state two immediate corollaries of Theorem
\ref{HKOmax_princ_theo}.
\begin{coro}\label{HKOuniq} (uniqueness)
   If $(u_1^1,\ldots,u_n^1)$ and $(u_1^2,\ldots,u_n^2)$ are solutions
   of \pP\ with the same initial condition
   $(u_{0,1},\ldots,u_{0,n})\in \big(C(\overline \Omega))^n$, then for
   all $i\in\no,\ u_i^1=u_i^2$ .\hskip 5pt $\blacksquare$
\end{coro}\medskip

\begin{coro}\label{HKOnoneg} (positivity)
   If $(u_1,\ldots,u_n)$ is the solution of \pP\ with the nonnegative initial condition $(u_{0,1},\ldots,u_{0,n})\in \big(C(\overline \Omega)\big)^n$, then for all $i\in\no,$ $u_i\geqslant 0$. Moreover, for all $i\in\no$, such that $u_{0,i}\geqslant 0$ and $u_{0,i} \not\equiv 0$, $u_i > 0$ in $\Omega$. \hskip 5pt $\blacksquare$
\end{coro}

%------------------------------------------------------------------------------%
\section{Contraction property}\label{HKOSecContr}
%------------------------------------------------------------------------------%

The purpose of this section is to show a contraction in $\big(L^1(\Omega)\big)^n$ pro\-per\-ty for the solutions of \pP\ with the initial conditions belonging to $\big(L^\infty(\Omega)\big)^n$. The main steps of the proof rely upon arguments due to \cite{bertsch_hilhorst} and \cite{osher_ralston}.

We first introduce some notation. We suppose that the functions $(u_1^1,\ldots,u_n^1)$ and
$(u_1^2,\ldots,u_n^2)$ are the solutions of \pP\  with the initial conditions $(u_{0,1}^1,\ldots,u_{0,n}^1)$
and $(u_{0,1}^2,\ldots,u_{0,n}^2)$, respectively. Define
   \begin{equation}\label{HKOmess2}
      (U_1,\ldots,U_n) := (u_1^1-u_1^2,\ldots,u_n^1-u_n^2).
   \end{equation}
Then
\begin{equation}\label{HKOprobl_U}
\begin{split}
      &(U_i)_t = \dvr \big(\sigma_i \nabla U_i + U_i \nabla \psi_i \big) \\
         &+\ \alpha_i\sum_{j=1}^n\lambda_{ij} \big(r_j(u_j^1(x,t),x) -r_j(u_j^2(x,t),x)\big) \quad\text{in}\quad Q_T,\\
            &\sigma_i\frac{\prt U_i}{\prt\nu}+U_i \frac{\prt \psi_i}{\prt\nu} = 0 \quad\text{on} \quad \partial\Omega\times (0,T),\\
      &U_i(x,0)=  U_{0,i}(x) \quad \text{for}\quad x\in \Omega,
   \end{split}
   \end{equation}
together with
\begin{equation}\label{HKOmess1}
U_{0,i}=u_{0,i}^1-u_{0,i}^2,
\end{equation}
for each $i\in\no$.\\

Next we prove the following contraction in $L^1$ property.

\begin{theo}\label{HKOWeakContrTh}
   For  all $t>0$,
   \begin{multline} \label{HKOWeakContr}
      \frac{1}{\alpha_1}\|U_1(\cdot,t)\|_{L^1(\Omega)} + \ldots +
      \frac{1}{\alpha_n}\|U_n(\cdot,t)\|_{L^1(\Omega)} \\
      \leqslant\frac{1}{\alpha_1}\|U_{0,1}(\cdot)\|_{L^1(\Omega)} + \ldots +
      \frac{1}{\alpha_n}\|U_{0,n}(\cdot)\|_{L^1(\Omega)},
   \end{multline}
   where $U_i$ and $U_{0,i},\ i\in \no$, are defined  by \eqref{HKOmess2} and \eqref{HKOmess1}, respectively.
\end{theo}
{\hskip -\parindent \bf Proof \hskip 5pt}
Dividing each partial differential equation of \eqref{HKOprobl_U} by $\alpha_i$ and summing them up, we obtain
   \begin{equation*}
   \begin{split}
      \frac{\text d}{\text dt}\bigg( \sum_{i=1}^{n}\frac{1}{\alpha_i} & U_i \bigg) =  \sum_{i=1}^n \frac{1}{\alpha_i}\,\dvr\left(\sigma_i\nabla U_i+U_i\nabla \psi_i\right) \\
      & +\sum_{i=1}^n \sum_{j=1}^n \lambda_{ij} \Big(r_j(u_j^1(x,t),x) - r_j(u_j^2(x,t),x) \Big) \\
      = & \sum_{i=1}^n \frac{1}{\alpha_i}\,\dvr\left(\sigma_i\nabla U_i+U_i\nabla \psi_i\right) \\
      & +\sum_{j=1}^n \bigg \{ \Big(r_j(u_j^1(x,t),x) - r_j(u_j^2(x,t),x) \Big)\sum_{i=1}^n \lambda_{ij} \bigg \} \\
      = & \sum_{i=1}^n \frac{1}{\alpha_i}\,\dvr\left(\sigma_i\nabla U_i+U_i\nabla \psi_i\right),
   \end{split}
   \end{equation*}
where we have used Hypothesis 2.

This, together with the boundary conditions \eqref{HKOmain_bd_con}, implies the conservation in time of the quantity
   \begin{equation}\label{HKOmess3}
     \frac{\text d}{\text dt} \sum_{i=1}^n \frac{1}{\alpha_i} \int_\Omega U_i(x,t) \dx =0.
   \end{equation}

Let us look closer at the nonlinear term
in \eqref{HKOprobl_U}. We can write, for fixed  index $i$
   \begin{multline*}
      \sum_{j=1}^n \lambda_{ij} \big( r_j(u_j^1(x,t),x)- r_j(u_j^2(x,t),x)\big)\\
         =\sum_{j=1}^n \lambda_{ij} U_j\int_0^1 \frac{\partial}{\partial u}r_j(\theta u_j^1 + (1-\theta)u_j^2,x)\text d\theta = \sum_{j=1}^n A_{ij} U_j.
   \end{multline*}

 Freezing the functions $u_i^k$ for $i\in\no,\ k\in\{1,2\}$, we deduce that the functions $U_1,\ldots,U_n$ satisfy a system of the form
   \begin{equation}\label{HKOlinear_eq}
      (U_i)_t =  \dvr \Big(\sigma_i \nabla U_i + U_i \nabla \psi_i \Big) + \sum_{j=1}^n A_{ij} U_{j} \quad \text{in}\quad Q_T,
   \end{equation}
with the boundary and initial conditions
   \begin{equation}\label{HKOlinear_cond}
   \begin{split}
      &\sigma_i\frac{\prt U_i}{\prt\nu}+U_i \frac{\prt \psi_i}{\prt\nu} = 0\quad\text{on}\quad \partial\Omega\times (0,T),\\
      &U_i(x,0)= U_{0,i}(x), \quad x \in \Omega.
   \end{split}
   \end{equation}
for $i\in\no$, where $A_{ij}$ are functions of space and time.

In order to make the notation more concise, we write
   \begin{equation*}
   \begin{split}
      \vec U_0 =        &\,\big(U_{0,1},\ldots,U_{0,n}\big),\\
      \vec U =          &\,\big(U_1,\ldots,U_n \big),\\
      \vec U_0^{\pm} =  &\,\big(U_{0,1}^{\pm},\ldots,U_{0,n}^{\pm}\big),\\
      \vec U^{\pm} =    &\,\big(U_1^{\pm},\ldots,U_n^{\pm}\big),\\
   \end{split}
   \end{equation*}
where $s^+=\max\{s,0\},\ s^-=\max\{-s,0\}$. By \eqref{HKOlinear_eq},\ \eqref{HKOlinear_cond}\ and Corollary \ref{HKOuniq} we can write $\vec U$ in the form
   \begin{equation*}
      (\vec U)(x,t)= \heat{t}\vec U_0(x) = \big( {\mathcal S}_1(t)\vec U_0,
      \ldots,{\mathcal S}_n(t)\vec U_0 \big)(x)
   \end{equation*}
with some operator ${\mathcal S}(t)$.
We set
\begin{equation*}
   \big( W_1,\ldots,W_n \big) = -\big(U_1 \ee^{\nicefrac{\psi_1(x)}{\sigma_1}},\ldots,U_n \ee^{\nicefrac{\psi_n(x)}{\sigma_n}}\big),
\end{equation*}
and $\widetilde A_{ij} = A_{ij}\, \ee^{\nicefrac{\psi_i(x)}{\sigma_i}} \, \ee^{-\nicefrac{\psi_j(x)}{\sigma_j}}$. Then, the system of equations \eqref{HKOlinear_eq} can be expressed in the form
   \begin{equation}\label{HKOWEq}
      \big( W_i \big)_t = {\sigma_i}\ee^{\nicefrac{\psi_i(x)}{\sigma_i}}\dvr \Big( \ee^{-\nicefrac{\psi_i(x)}{\sigma_i}} \nabla W_i \Big)
      + \sum_{j=1}^n \widetilde A_{ij} W_{j} \quad \text{in}\quad Q_T,
   \end{equation}
with the boundary and initial conditions
   \begin{align}
      \label{HKObdW} \frac{\prt W_i}{\prt \nu} &= 0 \quad \text{on}\quad \partial \Omega \times (0,T),\\
      \label{HKOWinit} W_i(x,0) &= - U_{0,i} \ee^{\nicefrac{\psi_i(x)}{\sigma_i}}, \quad x \in \Omega,
   \end{align}
for $i\in\no$.

Next we show that the solutions $W_i$ of the problem \eqref{HKOWEq}\,--\,\eqref{HKOWinit} with nonpositive initial conditions are nonpositive in $\overline \Omega$ for all $t\in(0,T)$. To that purpose we consider the auxiliary problem \begin{gather}
   \label{HKOaux1} \big( W_i \big)_t - \vartheta_i(x) \dvr \Big(\zeta_i(x) \nabla W_i \Big) - \sum_{j=1}^n \gamma_{ij} W_{j} \leqslant 0 \quad \text{in}\quad Q_T,\\
   \label{HKOaux2} \frac{\prt W_i}{\prt \nu} \leqslant 0 \quad \text{on}\quad \partial \Omega \times (0,T),\\
   \label{HKOaux3} W_i(x,0) = W_{0,i}(x) \leqslant 0 \quad x \in \Omega,
\end{gather}
for $i\in\no$. We assume that $\vartheta_i(x)$ and $\zeta_i(x)$ are nonnegative in $\overline \Omega$ and that the coefficients $\gamma_{ij}$ satisfy the same assumptions as the coefficients $\lambda_{ij}$ in \Pp.
The following result holds.
\begin{lemm}\label{HKOsummarize}
   Let $(W_1,\ldots,W_n)$ be a smooth and bounded solution of the problem \eqref{HKOaux1}\,--\,\eqref{HKOaux3}  with nonpositive initial conditions $W_{0,i}$ on a time interval $[0,T]$. Then $W_i(x,t) \leqslant 0$ in $\overline \Omega \times (0,T]$. Moreover, for each $i\in\no$ such that $W_{0,i} \leqslant 0$ and $W_{0,i} \not \equiv 0,\ W_i < 0$ in $\overline \Omega \times (0,T]$.
\end{lemm}
{\hskip -\parindent \bf Proof \hskip 5pt}
   The result of Lemma \ref{HKOsummarize} follows from the fact that the system (\ref{HKOaux1}), (\ref{HKOaux2}), (\ref{HKOaux3}), with the inequalities $\{\leqslant\}$  replaced by the equalities $\{=\}$, is a cooperative system. However, for the sake of
   completeness, we present a proof below.
   We first remark that, in view of \cite[Remark (i), p. 191]{protter_weinberger}, one can always satisfy the condition
 \begin{equation}\label{HKOwar_h_2}
 \sum_{j=1}^n \gamma_{ij} \leqslant 0 \text{\ for all\ } i \in \no,
\end{equation}
   for the matrix of coefficients $\big(\gamma_{ij}\big)_{i,j=1}^n$ by performing the change of variables $\overline W_i = W_i \ee^{-ct}$ for all $i\in \no$ and $c > 0$ large enough.\\
   %From the regularity of each $W_i$ there exists $M > 0$ such that $W_i - M \leqslant 0$ in $\overline \Omega \times [0,T]$ for all $i \in \no$. Then
   Thanks to the regularity of each $W_i$, we can apply Theorem 15, p.\ 191 from {\cite{protter_weinberger}} to conclude that $W_i - M \leqslant 0$ in $\overline \Omega \times [0,T]$ for some $M\> 0$ and all $i \in \no$. In fact, we can deduce that $W_i - M < 0$ in $\overline \Omega\times(0,T)$.\\
   Indeed, if for some $k \in \no$, $W_k = M$ in an interior point $(\tilde x, \tilde t)\in \Omega \times (0,T)$, then Theorem 15, p.\ 191 in {\cite{protter_weinberger}} implies that $W_k \equiv M$ for all $0 \leqslant t < \tilde t$, which is impossible since $W_k(x,0) \leqslant 0$.
   If the maximum $M$ of $W_k$ is attained at a boundary point $P\in \prt \Omega \times (0,T)$ then either there exists an open ball $K \subset \Omega \times (0,T)$ such that $P\in \prt K$ and $W_k - M < 0$ in $K$, and the last part of Theorem 15, p.\ 191 in {\cite{protter_weinberger}} contradicts the boundary inequality \eqref{HKOaux2}, or for all open balls $K \subset \Omega \times (0,T)$ such that $P\in \prt K$ there exists a point $(\tilde x, \tilde t) \in K$ such that $W_i(\tilde x, \tilde t) = M$, and we proceed as in the case before. \\
   Hence, there exists $\widetilde M > 0$, such that $W_i \leqslant \widetilde M < M$ in
   $\overline \Omega \times [0,T]$ for all $i \in \no$. Then we can repeat the reasoning for all $M > 0$
    until $M = 0$. Indeed, if this would not be the case, we  find the least real number $\overline M > 0$,
    with $W_i \leqslant \overline M \leqslant \widetilde M$ in $\overline \Omega \times [0,T]$, which leads again to the existence of a real number $0 \leqslant \widehat M < \overline M$ with the same property.
   This contradicts the fact that $\overline M$ was defined as the least such real number. \hskip 5pt $\blacksquare$ \medskip

Since the functions $u_i^1,\ u_i^2$ are bounded on $\overline \Omega \times [0,T]$, it follows that the functions $W_i$ are bounded on $\overline \Omega \times [0,T]$ for all $i\in \no$.\\
Then we are in a position to apply Lemma \ref{HKOsummarize} with $\vartheta_i(x) = \ee^{\nicefrac{\psi_i}{\sigma_i}}$, $\zeta_i(x) = \sigma_i \ee^{-\nicefrac{\psi_i}{\sigma_i}}$ and $\gamma_{ij} = \widetilde A_{ij}$ for $i,j \in \no$. We deduce that the solutions $W_i$ of the problem \eqref{HKOWEq}\,--\,\eqref{HKOWinit} with nonpositive initial conditions are nonpositive in $\overline \Omega$ for all $t\in(0,T)$.

Next we remark that the above reasoning can be applied either with $\vec U_0$ replaced by $U^+_0$ or with $\vec U_0$ replaced by $U^-_0$. This permits to show that ${\mathcal S}_i(t)\vec U_0^+,{\mathcal S}_i(t)\vec U_0^-\geqslant 0$ and that
\begin{equation}\label{HKOformm}
   {\mathcal S}_i(t)\vec U_0^\pm > 0 \quad \text{if} \quad \vec U_0^\pm \not\equiv 0.
\end{equation}

We easily compute
   \begin{equation}\label{HKOnazwac01}
   \begin{split}
      \sum_{i=1}^n \frac{1}{\alpha_i} &\big \|  U_i(\cdot,t) \big \|_{L^1(\Omega)} - \sum_{i=1}^n \frac{1}{\alpha_i} \big \| U_{0,i}(\cdot) \big \|_{L^1(\Omega)} \\
      =&\ \sum_{i=1}^n \frac{1}{\alpha_i}\big \| {\mathcal S}_i(t)\vec U_0^+-{\mathcal S}_i(t)\vec U_0^-\big \|_{L^1(\Omega)} -\sum_{i=1}^n \frac{1}{\alpha_i} \big \| U_{0,i}(\cdot)\big \|_{L^1(\Omega)}\\
      =&\ \sum_{i=1}^n \int_\Omega \frac{1}{\alpha_i}\Big \{ \max \big \{ {\mathcal S}_i(t) \vec U_0^+,{\mathcal S}_i(t)\vec U_0^-\big \}
   \end{split}
   \end{equation}
   \begin{equation}%\tag{\ref{HKOnazwac01}}
   \begin{split}
      &- \frac{1}{\alpha_i}\min \big \{{\mathcal S}_i(t)\vec U_0^+ ,{\mathcal S}_i(t) \vec U_0^- \big \} \Big \}\dx - \sum_{i=1}^n \frac{1}{\alpha_i}\int_\Omega \big \{U_{i,0}^+ + U_{i,0}^- \big \}\dx \\
      =&\ \sum_{i=1}^n \int_\Omega \frac{1}{\alpha_i} \big( {\mathcal S}_i(t)\vec U_0^+ + {\mathcal S}_i(t)\vec U_0^-\big)\dx - \sum_{i=1}^n \frac{1}{\alpha_i}\int_\Omega \big \{U_{i,0}^++U_{i,0}^-\big \}\dx \\
      & - 2 \sum_{i=1}^n\int_\Omega \frac{1}{\alpha_i} \min \big \{{\mathcal S}_i(t)\vec U_0^+,{\mathcal S}_i(t)\vec U_0^-\big \}\dx \\
      =&\ - 2 \sum_{i=1}^n\int_\Omega \frac{1}{\alpha_i}\min \big \{{\mathcal S}_i(t)\vec U_0^+,{\mathcal S}_i(t)\vec U_0^-\big \}\dx \leqslant 0,
   \end{split}
   \end{equation}
which completes the proof of (\ref{HKOWeakContr}).\hskip 5pt $\blacksquare$ \medskip

\begin{coro}\label{HKOStrContrTh}
   Let $(u_{0,1}^1, \ldots, u_{0,n}^1),\ (u_{0,1}^2, \ldots, u_{0,n}^2) \in \big(C(\overline \Omega)\big)^n$
   be as in Theorem \ref{HKOWeakContrTh}. Moreover, let us assume that for at least one index $k \in \no$
   the difference $u^1_{0,k} - u^2_{0,k}$ changes the sign. Then, the inequality \eqref{HKOWeakContr} is strict for all $t > 0$, so that %the semigroup
   solution satisfies a strict contraction property.
\end{coro}

%------------------------------------------------------------------------------%
\section{Large time behavior of solutions}\label{HKOlargetime}
%------------------------------------------------------------------------------%

In this section we assume the existence and uniqueness of a positive solution $\vec v = (v_1, \ldots, v_n) \in \big( C(\overline\Omega) \cap C^2(\Omega)\big)^n$ of the elliptic problem
   \begin{align}
      \dvr \big(\sigma_i\nabla v_i + v_i \nabla \psi_i \big) + \alpha_i \bigg( \sum_{j=1}^n \lambda_{ij} r_j \big(v_j(x),x \big) \bigg) &= 0 \quad \text{in}\quad \Omega, \label{HKOStacEq}\\
      \sigma_i \frac{\prt v_i}{\prt \nu} + v_i \frac{\prt \psi_i}{\prt \nu} &= 0 \quad \text{on}\quad \prt \Omega, \label{HKOStacInit}\\
      \sum_{i=1}^n \frac{1}{\alpha_i} \int_\Omega v_i(x) \dx &= 1, \label{HKOWarCalk}
   \end{align}
for $i \in \no$.
\medskip

\begin{defi}
   We say that a vector function $\vec v = (v_1,\ldots,v_n) \in \big( C(\overline\Omega )\big)^n$ is nonnegative (resp.\ positive) if $v_i(x) \geqslant 0$ (resp.\ $v_i(x) > 0$) for all $x\in \overline \Omega$ and all $i\in\no$.
\end{defi}
Next we introduce the semigroup notation for the unique solution of Pro\-blem~(P), namely
\begin{equation*}
   \vec u(t) = \heatT{t} \vec u_0 = \Big( \mathcal T_1(t)\,\vec u_0,
      \ldots,\mathcal T_n(t)\,\vec u_0 \Big),
\end{equation*}
with the initial data $\vec u_0 \in \big(C(\overline\Omega)\big)^n$. The method of the proof is based upon
an idea of Osher and Ralston \cite{osher_ralston}. It mainly exploits the contraction properties for the
nonlinear semigroup $\heatT{t}$ given by Theorem \ref{HKOWeakContrTh} and Corollary \ref{HKOStrContrTh}. A similar reasoning was developed in other contexts\ by Bertsch and Hilhorst \cite{bertsch_hilhorst}, Hilhorst and Hulshof \cite{hilhorst_hulshof} and Hilhorst and Peletier \cite{hilhorst_peletier}.
\medskip

\noindent
We suppose there exists a set $\mathscr H \subset \big( C(\overline\Omega) \cap C^2(\Omega) \big)^n$ of positive stationary solutions with the following property which we denote by $\mathscr S$:\\
\textit{For each $\vec f = (f_{1}, \ldots, f_{n}) \in \big( C(\overline\Omega) \cap C^2(\Omega) \big)^n$ either $\vec f \in \mathscr H$ or there exists $(\xi_{1}, \ldots, \xi_{n}) \in \mathscr H$, such that $f_i-\xi_i$ changes the sign for at least one index $i\in\no$.}
\medskip

One can prove that a set $\mathscr H$ satisfying Property $\mathscr S$ exists in at least two cases:
\begin{itemize}
   \item[\textit{i)}]
      In the case of the system \eqref{HKOProblemBothe} where the Robin boundary conditions reduce to the homogeneous Neumann boundary conditions, the set $\mathscr H$ is given by
      \begin{multline*}
         \rule{20pt}{0pt} \mathscr H = \Big \{ (a,b):\ a>0,\ b= r_B^{-1}(r_A(a))\\
         \text{and}\ \frac{a}{\alpha} + \frac{b}{\beta} = \int_\Omega \Big(\frac{u}{\alpha} + \frac{v}{\beta} \Big) \dx \Big \}.\rule{20pt}{0pt}
      \end{multline*}
      For more details we refer to \cite{bothe_hilhorst}.
   \item[\textit{ii)}]
      In the case of the molecular motor with a linear $n$-component system the set $\mathscr H$ is given by
      \begin{equation*}
         \mathscr H = \big \{ c \vec v:\ c\in \mathbb R^+ \big\},
      \end{equation*}
      where $\vec v$ is a unique solution of the elliptic problem \eqref{HKOStacEq}\,--\,\eqref{HKOWarCalk}.
\end{itemize}

\begin{prop}\label{HKOcontin}
   The continuum $\mathscr H$ is such that for each
   \begin{equation*}
      \vec f = (f_{1}, \ldots, f_{n}) \in \big( C(\overline\Omega) \cap C^2(\Omega) \big)^n
   \end{equation*}
   either $\vec f \in \mathscr H$, or there exists $(\xi_{1}, \ldots, \xi_{n}) \in \mathscr H$ such that $f_i-\xi_i$ changes the sign for at least one index $i\in\no$.
\end{prop}
{\hskip -\parindent \bf Proof \hskip 5pt}
\begin{itemize}
\item[\textit{i)}]
   In the case of system \eqref{HKOProblemBothe} the proof is rather obvious since the continuum $\mathscr H$ is composed of constant pairs.
\item[\textit{ii)}]
   In the case of the molecular motor, let us assume that $\vec f \not\in \mathscr H$. Then there does not exist any positive constant $c$ such that $c\,\vec v = \vec f$. In particular, there exists an index $i\in\no$ such that $v_i$ is not proportional to $f_i$, or in other words $c v_i \neq f_i$ for all $c>0$. Without loss of generality we can assume that the first coordinate has this property. Let $x_0 \in \Omega$ be arbitrary. Since $v_1$ is strictly positive in $\overline \Omega$, we can define
   \begin{equation*}
      c_0 = \frac{f_1(x_0)}{v_1(x_0)},
   \end{equation*}
   so that
   \begin{equation*}
      \big(f_1 - c_0\,v_1 \big)(x_0) = 0.
   \end{equation*}
   Let $\mathcal Z = \big \{ x\in\overline\Omega:\ \big(f_1 - c_0\, v_1\big)(x) = 0 \big \}$. From the continuity of $f_1$ and $v_1$, $\mathcal Z$ is closed as a subset of $\Omega$. If there exist $x_1,\, x_2 \in \mathcal Z^c$, such that $\big(f_1 - c_0\,v_1 \big)(x_1)$ and $\big(f_1 - c_0 \,v_1 \big)(x_2)$ are of different signs, then the proof is complete. Now suppose that $\big( f_1 - c_0\,v_1 \big)(x)$ is positive for all $x \in \mathcal Z^c$. In particular
   \begin{equation*}
      \big(f_1 - c_0 \, v_1 \big)(\tilde x) = d > 0
   \end{equation*}
   for some fixed $\tilde x \in \mathcal Z^c$. Then choosing $\displaystyle \eps = \frac{d}{2 v_1(\tilde x)}$ we see that
   \begin{equation*}
      \big(f_1 - (c_0 + \eps ) v_1 \big)(\tilde x) = \frac{d}{2} > 0.
   \end{equation*}
   However
   \begin{equation*}
      \big(f_1 - (c_0 + \eps ) v_1 \big)(x_0) < 0.
   \end{equation*}
   We proceed similarly when $\big( f_1 - c_0\,v_1 \big)(x)$ is negative for all $x \in \mathcal Z^c$.
   \hskip 5pt $\blacksquare$ \medskip
\end{itemize}

In the sequel we suppose that the initial data $\vec u_0 = (u_{0,1},\ldots,u_{0,n}) $ from $\big( C(\overline \Omega)\big)^n$ also satisfy the following property:
\begin{equation}\label{HKOcokolwiek}
   \text{\textit{There exists\ }} \vec h \in \mathscr H \text{\textit{\ such that\ }} 0 \leqslant \vec u_0 \leqslant \vec h \quad \text{in} \quad \overline \Omega,
\end{equation}
and remark that this property is satisfied in both the cases \textit{(i)} and \textit{(ii)}.

\begin{prop}\label{HKOnadrozwiazanie}
   Let $\vec u_0 = (u_{0,1},\ldots,u_{0,n}) \in \big( C(\overline \Omega)\big)^n$ satisfy the property \eqref{HKOcokolwiek}. Then the solution $(u_1,\ldots,u_n)$ of \Pp~ is such that $0 \leqslant \vec u(t) \leqslant \vec h$ for all $t > 0$.
\end{prop}
{\hskip -\parindent \bf Proof \hskip 5pt}
   We remark that $\vec 0$ is a subsolution of \Pp~ and that $\vec h$ is a supersolution, and apply Theorem \ref{HKOmax_princ_theo}.
\hskip 5pt $\blacksquare$ \medskip

Next we prove the main result of this section. To that purpose we first define the norm $\pmb{\big \|} \cdot \pmb{\big \|_1}$ by
\begin{equation*}
   \pmb{\big \|} \vec f\ \pmb{\big \|_1} = \sum_{i=1}^n \frac{1}{\alpha_i} \big \| f_i \big \|_{L^1(\Omega)}\hskip 5pt.
\end{equation*}
Note that this norm is equivalent to the usual product norm in the space $\big( L^1(\Omega)\big)^n$.
\begin{theo}\label{HKOAsympt}
   For all nonnegative $\vec u_0 = (u_{0,1},\ldots,u_{0,n})\in \big( C(\overline \Omega)\big)^n$ there exists $\vec f = (f_{1}, \ldots, f_{n}) \in \mathscr H$, such that
   \begin{equation*}
      \lim_{t \rightarrow \infty} \pmb{\big \|}\ \heatT{t} \vec u - \vec f\ \pmb{\big \|_1} = 0.
   \end{equation*}
\end{theo}
{\hskip -\parindent \bf Proof \hskip 5pt}\\
   The proof consists of several steps. To begin with we define the $\omega$-limit set
   \begin{multline}\label{HKOomegaset}
      \omega(\vec u_0) = \Big \{ \vec g \in \big( L^1(\Omega)\big)^n: \ \text{there exists a sequence \ } t_k \rightarrow \infty
      \\ \text{\ as\ } k \rightarrow \infty,\text{\ such that\ }\lim_{k \rightarrow \infty} \pmb{\big \|} \heatT{t_k}\vec u_0 - \vec g\ \pmb{\big \|_1} = 0 \Big \},
   \end{multline}
The organization of the proof is as follows. First we show that $\omega(\vec u_0)$ is not empty. In the second step we define the Lyapunov functional
   \begin{equation*}
      \calV(\vec \xi) = \pmb{\big \|} \vec \xi - \vec w \pmb{\big \|_1},
   \end{equation*}
   where $\vec w$ is a stationary solution and check that it is constant on $\omega(\vec u_0)$. We then deduce that $\omega(\vec u_0) \subset \mathscr H$, and finally prove that $\omega(\vec u_0)$ consists of exactly one function.

   \vskip 10pt
   \noindent \textbf{Step 1.} {\textit{$\omega(\vec u_0)$ is not empty}}.\\
   Let $\eps >0$ be arbitrary. Suppose that $\Omega' \subset\subset \Omega$ satisfy
   \begin{equation*}
      \big | \Omega \setminus \Omega' \big | \leqslant \frac{\eps}{2K}.
   \end{equation*}
   and set
   \begin{equation}\label{HKOK}
      K = \sum_{i=1}^n \frac{2}{\alpha_i} \|h_i\|_{C(\overline \Omega)},
   \end{equation}
   where $\vec h$ has been introduced in \eqref{HKOcokolwiek}.
   We have already proved in Proposition \ref{HKOnadrozwiazanie} that $\heatT{t} \vec u_0$ is bounded in $\big( L^\infty(\Omega)\big)^n$. Therefore there exist a vector function $\vec g \in \big{(L^\infty}(\Omega))^n$ and a sequence $\{\vec u(t_k)\}$ such that
   \begin{equation}\label{HKOdiffomega}
      \vec u(t_k) \rightharpoonup \vec g \mbox{~~weakly in ~~} (L^2(\Omega))^n,
   \end{equation}
   as $t_k \rightarrow \infty$. Next we deduce from \cite[Chap.\ \texttt{III},\ Theorem 10.1]{lady_solo_ural} that there exists a positive constant $C$ such that
   \begin{equation*}
      \big| u_i(x_1,t) - u_i(x_2,t) \big| \leqslant C|x_1 - x_2|^{\alpha}
   \end{equation*}
   for all $x_1,\ x_2 \in \Omega'$ and all $t>0$. Therefore, it follows from the Ascoli-Arzel\`a Theorem (see, e.g., \cite[Theorem 1.33]{adams}) that $\vec u(t_k) \rightarrow \vec g$ as $t_k\rightarrow \infty$, uniformly in $\overline \Omega'$. We choose $t_0$ large enough such that for all $t_k \geqslant t_0$
   \begin{equation}\label{HKOnormdiff}
      \pmb{\big \|} \vec u(\cdot,t_k) - \vec g(\cdot) \pmb{\big \|}_{\pmb{1},\Omega'} \leqslant \frac{\eps}{2},
   \end{equation}
   where $\pmb{\big \|} \cdot \pmb{\big \|}_{\pmb{1},\Omega'}$ corresponds to the $L^1$ norm in $\Omega'$. We deduce that, in view of \eqref{HKOK} and \eqref{HKOdiffomega} that
   \begin{equation*}
      \pmb{\big \|} \vec u(\cdot,t_k) - \vec g(\cdot) \pmb{\big \|}_{\pmb{1},\Omega \setminus \Omega'} \leqslant K \big| \Omega \setminus \Omega' \big| \leqslant \frac{\eps}{2},
   \end{equation*}
   which together with \eqref{HKOnormdiff} yields
   \begin{equation*}
      \pmb{\big \|} \vec u(\cdot,t_k) - \vec g(\cdot) \pmb{\big \|_1} \leqslant \eps.
   \end{equation*}

   \vskip 10pt
   \noindent \textbf{Step 2.} \textit{$\omega(\vec u_0) \subset \mathscr H$.}\\
   Indeed, let $\vec g \in \omega(\vec u_0)$ and suppose $\vec g \notin \mathscr H$. According to Proposition \ref{HKOcontin} we can find a steady state solution $\vec w\in \mathscr H$, such that at least one component of $\vec w - \vec g$ changes the sign. Without loss of generality we can assume that it happens for the first component, namely that\ $f_1 - w_1$ changes the sign. We remark that, by the contraction property in Theorem \ref{HKOWeakContrTh}, the  functional
   \begin{equation*}
      \calV(\vec \xi) = \pmb{\big \|} \vec \xi - \vec w \pmb{\big \|_1}
   \end{equation*}
   is a Lyapunov functional for \Pp, where $\vec \xi\in\big( L^1(\Omega)\big)^n$. Next we describe some of its properties.

   \vskip 10pt
   \noindent \textbf{Property (a)} {\textit{ The functional $\calV$ is constant on $\omega(\vec u_0)$}}.\\
   Since $\heatT{t} \vec w = \vec w$ and $\heatT{t}$ has the contraction property \eqref{HKOWeakContr}, the functional $\calV$ is nonincreasing in time along the trajectory $\heatT{t} \vec u_0$, which yields
   \begin{equation*}
   \begin{split}
      \calV \big( \heatT{t} \vec u_0 \big ) &= \pmb{\big \|} \heatT{t} \vec u_0 - \vec w \pmb{\big \|_1} \\
      &= \pmb{\big \|} \heatT{t} \vec u_0 - \heatT{t} \vec w \pmb{\big \|_1} \leqslant \pmb{\big \|} \vec u_0 - \vec w \pmb{\big \|_1} < \infty \hskip 5pt.
   \end{split}
   \end{equation*}
   Thus there exists a finite limit $\calV^\ast$ of $\calV \big( \heatT{t} \vec u_0 \big )$ as $t \rightarrow \infty$. Let $\vec h_1,\ \vec h_2 \in \omega(\vec u_0)$. We can find a sequence $t_k \rightarrow \infty$ as $k \rightarrow \infty$, such that
   \begin{equation*}
   \pmb{\big \|} \heatT{t_{2k}} \vec u_0 - \vec h_1 \pmb{\big \|_1} \rightarrow 0 \quad \text{and} \quad \pmb{\big \|} \heatT{t_{2k+1}} \vec u_0 - \vec h_2 \pmb{\big \|_1} \rightarrow 0,
   \end{equation*}
   as $k$ tends to $\infty$. It follows that $\calV \big( \vec h_1 \big) = \calV \big( \vec h_2 \big) = \calV^\ast$.

   \vskip 10pt
   \noindent \textbf{Property (b)} {\textit{The $\omega$-limit set $\omega(\vec u_0)$ is invariant with respect to the semigroup $\heatT{t}$}, namely if $\vec h\in\omega(\vec u_0)$, then for all $t>0$ also $\heatT{t} \vec h \in\omega(\vec u_0)$}.\\
   Let the sequence $t_k \rightarrow \infty$ as $k \rightarrow \infty$ be such that $\pmb{\big \|} \heatT{t_k} \vec u_0 - \vec h \pmb{\big \|_1} \rightarrow 0$. From the contraction property \eqref{HKOWeakContr}
   \begin{gather*}
      \pmb{\big \|} \heatT{t_k + t} \vec u_0 - \heatT{t} \vec h \pmb{\big \|_1} = \pmb{\big \|} \heatT{t}\heatT{t_k} \vec u_0 - \heatT{t} \vec h \pmb{\big \|_1} \\
      \leqslant \pmb{\big \|} \heatT{t_k} \vec u_0 - \vec h \pmb{\big \|_1}.
   \end{gather*}
   Since the last term above tends to $0$ as $k$ tends to $\infty$ this shows that $\heatT{t} \vec h \in\omega(\vec u_0)$.

   \vskip 10pt
   Now, remember that $\vec g\in\omega(\vec u_0)$ is such that $\vec g\notin \mathscr H$ and $\vec w \in \mathscr H$ is such that the first component of $\vec w - \vec g$ changes the sign in $\Omega$. Then, Corollary \ref{HKOStrContrTh} yields
   \begin{multline*}
      \rule{20pt}{0pt}\calV(\heatT{t} \vec g) = \pmb{\big \|} \heatT{t} \vec g - \vec w \pmb{\big \|_1} \\
      = \pmb{\big \|} \heatT{t} \vec g - \heatT{t} \vec w \pmb{\big \|_1} < \pmb \|\vec g - \vec w \pmb{\|_1} = \calV(\vec g),\rule{20pt}{0pt}
   \end{multline*}
   for all $t>0$, which contradicts Property (a). Therefore $\vec g \in \mathscr H$.

   \vskip 10pt
   \noindent \textbf{Step 3.} \textit{The set $\omega(\vec u_0)$  contains only one element}.\\
   Suppose that $\vec g_1,\vec g_2 \in \omega(\vec u_0)$. Then we can find two sequences $t_k,\,s_k$ tending to $\infty$ as $k \rightarrow \infty$, such that $s_k \leqslant t_k$ and $\pmb{\big \|} \heatT{t_k} \vec u_0 - \vec g_1 \pmb{\big \|_1},\, \pmb{\big \|} \heatT{s_k} \vec u_0 - \vec g_2 \pmb{\big \|_1} \rightarrow 0$ as $t_k \rightarrow \infty$. Since $\omega(\vec u_0) \subset \mathscr H$, it follows that
   \begin{gather*}
      \pmb \|\vec g_1 - \vec g_2 \pmb{\|_1} \leqslant \pmb{\big \|} \heatT{t_k} \vec u_0 - \vec g_1 \pmb{\big \|_1} + \pmb{\big \|} \heatT{t_k} \vec u_0 - \vec g_2 \pmb{\big \|_1}  \\
      =\, \pmb{\big \|} \heatT{t_k} \vec u_0 - \vec g_1 \pmb{\big \|_1}%\\
      + \pmb{\big \|} \heatT{t_k - s_k} \heatT{s_k} \vec u_0 - \heatT{t_k - s_k} \vec g_2 \pmb{\big \|_1}\\
       \leqslant\pmb{\big \|} \heatT{t_k} \vec u_0 - \vec g_1 \pmb{\big \|_1} + \pmb{\big \|} \heatT{s_k} \vec u_0 - \vec g_2 \pmb{\big \|_1} \hskip 3pt,
   \end{gather*}
   which tends to $0$ as $k\rightarrow \infty$.\hskip 5pt $\blacksquare$

%------------------------------------------------------------------------------%
\section{Stationary so\-lu\-tions for the linear molecular motor problem}
%------------------------------------------------------------------------------%

In this section we show the existence and the uniqueness (up to a multiplicative constant) of the classical stationary solution of the problem for the molecular motor. We suppose that $\Omega$ is an open  bounded subset of $\mathbb R^d$ with smooth boundary $\prt \Omega$.

We consider the linear system
\begin{equation}\label{HKOnStacEq1}
    \dvr \big(\sigma_i \nabla v_i(x) + v_i(x) \nabla \psi_i(x) \big) + \sum_{j=1}^n a_{ij} v_j(x) = 0 \quad \text{in} \quad \Omega,
\end{equation}
where $i\in \no$, $n>1$.
\noindent The system \eqref{HKOnStacEq1} is supplemented with the Robin boundary conditions
\begin{equation}\label{HKOnStacBd1}
   \sigma_i\frac{\partial v_i}{\partial \nu} + v_i\frac{\partial \psi_i}{\partial \nu} = 0 \quad \text{on}\quad \partial \Omega,
\end{equation}
where $i \in \no$. Thus, the problem can be written as
\begin{equation*}
   {\cal A}\vec v = 0,
\end{equation*}
with a linear operator $\cal A$ in a suitable Banach space $\cal X$ of functions on $\Omega$, to be made precise later. Moreover, we impose the integral constraint
\begin{equation}\label{HKOintFirst}
   \sum_{i=1}^n \int_\Omega v_i(x) \dx = 1.
\end{equation}

The adjoint problem ${\cal A}^\ast\vec \f=0$ to \eqref{HKOnStacEq1}, in a dual space ${\cal X}^\ast$, is now
\begin{equation}\label{adj}
   \sigma_i \Delta \f_i - \nabla \psi_i \cdot \nabla \f_i + \sum_{j=1}^n a_{ji} \f_j = 0, \quad \text{in}\quad \Omega,
\end{equation}
with the Neumann boundary conditions for each $i=1,\,\dots,\,n$
\begin{equation}\label{Neu}
   \frac{\partial \f_i}{\partial\nu} = 0 \quad \text{on}\quad \partial \Omega.
\end{equation}
Since $\sum_{j=1}^na_{ji}=0$, the problem \eqref{adj} has the obvious solution
\begin{equation}\label{phi}
   \vec \f = (\f_1,\dots,\f_n) = (1,\dots,1).
\end{equation}
We are going to apply the Krein-Rutman theorem on the first eigenvalues and eigenvectors of positive operators, and this will permit us to conclude that the problem \eqref{HKOnStacEq1}--\eqref{HKOnStacBd1} has a one-dimensional space of solutions. Therefore, under the additional constraint \eqref{HKOintFirst}, the original problem \eqref{HKOnStacEq1}--\eqref{HKOnStacBd1} has a unique solution.\\
Perthame and Souganidis sketched this argument for $n>1$ and $d=1$ in \cite{ps}.

\begin{theo}\label{Th4}
   Under the assumption $\sum_{j=1}^na_{ji}=0$, there exists a unique smooth solution $\vec v$ of the system (\ref{HKOnStacEq1})--(\ref{HKOintFirst}).
\end{theo}
Before proving Theorem \ref{Th4} we recall some basic definitions as well as  the Krein-Rutman theorem from \cite[Ch. VIII, p. 188--191]{dl}.

\begin{defi}[Reproducing cone]
   We say that a closed set $K$ in $\cal X$ is a cone, if it possesses the following properties:
   \begin{itemize}
      \item[i)] $0 \in K$,
      \item[ii)] $u,\ v \in K \Longrightarrow \alpha u + \beta v \in K$, for all $\alpha,\ \beta \geqslant 0$,
      \item[iii)] $v \in K$ and $-v \in K \Longrightarrow v = 0$.
   \end{itemize}
   A cone $K \subset \cal X$ is said to be reproducing if ${\cal X} = K - K \equiv \big \{ k_1 - k_2:\ k_1,\ k_2 \in K \big \}$.
\end{defi}

\begin{defi}[Dual cone]
   If $K$ is a cone in $\cal X$, then the set $K^\ast \subset {\cal X}^\ast$ is said to be a dual cone if
   \begin{equation*}
      \langle f^\ast , v \rangle \geqslant 0,
   \end{equation*}
   for every $v \in K$.
\end{defi}

\begin{defi}[Strict positivity]\label{HKOstrPos}
   Let $\cal B$ be a linear operator on $\cal X$. Then $\cal B$ is said to be strongly positive if ${\cal B} v \in K^o$ for all $v \in K$ such that $v \neq 0$.
\end{defi}

\begin{theo}\label{KR}
   Let $K$ be a reproducing cone in a Banach space $\cal X$, with nonempty interior $K^o\neq\emptyset$, and let $\cal B$ be a strongly positive compact operator on $K$ in a sense of Definition \ref{HKOstrPos}. Then the spectral radius of $\cal B$, $r(\cal B)$, is a simple eigenvalue of $\cal B$ and $\cal B^\ast$, and their associated eigenvectors belong to $K^o$ and $(K^\ast )^o$. More precisely, there exists a unique associated eigenvector in $K^o$ (resp. $(K^\ast)^o$) of norm $1$. Furthermore, all other eigenvalues are strictly less in absolute value than $r(\cal B)$.
\end{theo}

{\hskip -\parindent \bf Proof \hskip 5pt}
   We will apply Theorem \ref{KR} to the space ${\cal X} = \big( C(\overline\Omega) \big)^n \subset \big(L^1(\Omega)\big)^n$ endowed with the usual supremum norm, and the operators
   \begin{gather*}
      {\cal B} = (\lambda I-{\cal A})^{-1}:\ {\cal X}\to{\cal X},\\
      {\cal B}^\ast = (\lambda I-{\cal A}^\ast)^{-1}:\ {\cal X}^\ast \to {\cal X}^\ast,
   \end{gather*}
   where $\lambda>0$ is a strictly positive real number to be fixed later.

   Let
   \begin{equation*}
      K = \big \{\vec u\in{\cal X}:\ u_i(x) \geqslant 0\ {\rm for\ each\ } x \in \overline \Omega,\ i=1,\dots,n \big \}.
   \end{equation*}
   We remark that $K$ is a reproducing cone, with nonempty interior
   \begin{equation*}
      K^o = \big \{ \vec u\in{\cal X}:\ \inf_{x \in \overline \Omega} u_i(x) > 0 ,\ i=1,\dots,n \big \}.
   \end{equation*}
   From the standard theory \cite[Theorem 2.1 and Theorem 3.1, Ch. 7]{lady_ural} for elliptic partial differential linear systems, the boundary value problem
   \begin{equation}\label{adj2}
      \sigma_i \Delta \f_i - \nabla \psi_i \cdot \nabla \f_i + \sum_{i=1}^n a_{ji} \f_j - \lambda \f_i = f_i\ \ {\rm in\ \ }\Omega,
   \end{equation}
   with the homogeneous Neumann conditions \eqref{Neu} on $\partial\Omega$, for $\lambda = \widetilde \lambda >0$ sufficiently large, has a solution $\vec \f = (\f_1,\dots,\f_n) \in {\cal X}$ for each $\vec f=(f_1,\dots,f_n)\in{\cal X}$. Moreover, if $f_i(x)\geqslant 0$ for each $i=1,\dots,n$, and $x\in\overline\Omega$, then $\f_i(x)\geqslant 0$ (in fact, $\f_i(x)>0$ in $\Omega$), which is a consequence of the maximum principle (cf. also Example 3 on p. 196--197 in \cite{dl}). Thus, the operator ${\cal B}^\ast = \big(\widetilde \lambda I-{\cal A}^\ast\big)^{-1}$ is a strongly positive and compact operator, and by Theorem \ref{KR}, the largest eigenvalue $\mu$ of $\cal B$ and ${\cal B}^\ast$ is simple.

   Since
   \begin{gather*}
      - \sigma_i \Delta \f_i + \nabla \psi_i \cdot \nabla \f_i - \sum_{j=1}^n a_{ji} \f_j + \widetilde \lambda \f_i = \widetilde \lambda \f_i\quad \text{in}\quad \Omega \\
      \frac{\partial \f_i }{\partial \nu} = 0 \quad \text{on}\quad \partial \Omega,
   \end{gather*}
   for all $i \in \no$, with $\vec \f = (\f_1,\ldots,\f_n) = (1,\ldots,1)$, and since $(1,\ldots,1) \in ( K^\ast )^o$, it follows that $\displaystyle \frac{1}{\widetilde \lambda} = r\Big( \big(\widetilde \lambda I - {\cal A}^\ast \big)^{-1} \Big)$ is a simple eigenvalue of the operator $\big(\widetilde \lambda I - {\cal A}^\ast \big)^{-1}$. Applying again Theorem \ref{KR}, we deduce that $\displaystyle \frac{1}{\widetilde \lambda}$ is the largest eigenvalue of the operator $\big(\widetilde \lambda I - {\cal A} \big)^{-1}$ and that it is simple, and that there exists $\vec v \in K^o \subset \cal X$ such that
   \begin{equation*}
      \Big(\widetilde \lambda I - {\cal A} \Big)^{-1} \vec v = \frac{1}{\widetilde \lambda}\, \vec v,
   \end{equation*}
   which is equivalent to
   \begin{equation*}
      {\cal A} \vec v = 0.
   \end{equation*}
   This proves the existence of the solution of the problem \eqref{HKOnStacEq1}--\eqref{HKOintFirst}.  \hskip 10pt $\blacksquare$

   %In fact, $\mu=1/\lambda$ since \eqref{phi} holds, and  $B'\vec\f=\mu\vec\f$ is equivalent to ${\cal A}'\vec\f=(\lambda-1/\mu)\vec\f$. What we really need is that $\mu$ is a simple eigenvalue of $B$, and $\vec v\in{\cal X}$ satisfying
   %\begin{equation*}
   %   B\vec v=\mu\vec v
   %\end{equation*}
   %belongs to $(K^o)'$. Since the above equation is equivalent to ${\cal A}\vec v=\vec v$, we conclude that a solution to \eqref{HKOnStacEq1}--\eqref{HKOnStacBd1} is unique up to  a multiplicative constant, so \eqref{HKOintFirst} is satisfied with a suitable normalization of $\vec v$.

%------------------------------------------------------------------------------%

%------------------------------------------------------------------------------%

\end{document}